\documentclass{article}
\usepackage{amsmath,amsthm,amsopn,amstext,amscd,amsfonts,amssymb,verbatim}
\usepackage{pdflscape}
\usepackage{natbib}

\def\diag{\mathop{\rm diag}\nolimits}
\def\tr{\mathop{\rm tr}\nolimits}

\def\build#1#2#3{\mathrel{\mathop{#1}\limits^{#2}_{#3}}}

\def\Vec{\mathop{\rm vec}\nolimits}
\def\Cov{\mathop{\rm Cov}\nolimits}
\def\Var{\mathop{\rm Var}\nolimits}

\def\E{\mathop{\rm \mathbf{E}}\nolimits}
\def\P{\mathop{\rm \mathbf{P}}\nolimits}

\renewenvironment{abstract}
                 {\vspace{6pt}
                  \begin{center}
                  \begin{minipage}{5in}
                  \centerline{\textbf{Abstract}}
                  \noindent\ignorespaces
                 }
                 {\end{minipage}\end{center}}

\theoremstyle{definition}

\title{\Large \textbf{Multiple response optimisation: Multiobjective stochastic programming methods}}
\author{
  \textbf{Jos\'e A. D\'{\i}az-Garc\'{\i}a} \thanks{Corresponding author\newline
    }\\
  {\normalsize Department of Statistics and Computation} \\
  {\normalsize 25350 Buenavista, Saltillo, Coahuila, Mexico} \\
  {\normalsize E-mail: jadiaz@uaaan.mx} \\[2ex]
  \textbf{Mahdi Bashiri} \\
  {\normalsize Industrial Engineering Department} \\
  {\normalsize Shahed University, Faculty of Engineering} \\
  {\normalsize P. O. Box 18155/159, Teheran, Iran}\\
  {\normalsize E-mail: bashiri@shahed.ac.ir}\\
}
\date{}
\begin{document}
\maketitle

\begin{abstract}
The multiresponse surface problem is modelled as one of multiobjective stochastic optimisation,
and diverse solutions are proposed. Several crucial  differences are highlighted between this
approach and others that have been proposed. Finally, in a numerical example, some particular
solutions are applied and described in detail.
\end{abstract}

\section{Introduction}\label{Sec1}

Many (perhaps most) real-world design problems are in fact multiobjective optimisation problems
in which the designer seeks to optimise simultaneously several performance attributes of a
design,  and an improvement in one objective is often only gained at the cost of deteriorations
in others,  and so a trade-off is necessary. Similar situations are met in the study  of
natural phenomena and in experimental trials.

Moreover,  the response variables,  perhaps the controllable variables and even some parameters
involved in these studies may have a random (or stochastic) character.

A very useful statistical tool in the study of these designs,  phenomena and experiments is
that of the response surfaces methodology,  in its multivariate version. This approach makes it
possible to determine an analytical relationship between the response and control variables,
through a process of continuous improvement and optimisation. Similarly,  it allows us to
obtain an approximate vector function (termed the multiresponse surface or predicted response
vector) with a smaller amount of data and fewer experimental runs,  see \citet{kc:87},
\citet{k:08a},  \citet{mm:02} and \citet{k:08b}.

Although the response variables are random and in consequence the estimated multiresponse
surface contains shape parameters that are estimated in the regression stage (i.e. they are
random),  the process was initially considered as one of deterministic optimisation,  see
\citet{b:75} among others. Subsequently,  this randomness or uncertainty was taken into account
in the multiobjective optimisation process in different ways,  and at different stages ,  see
\citet{kc:81},  \citet{kc:87} \citet{chh:01},  \citet{n:08},  \citet{aea:08} and
\citet{hbdn:10}.

There exists a well-grounded and documented theory -Stochastic Optimisation- that addresses the
following general problem:
\begin{equation}\label{eq0}
  \begin{array}{c}
  \build{\min}{}{\mathbf{x}}
  \left(
  \begin{array}{c}
    h_{1}(\mathbf{x}, \boldsymbol{\xi}) \\
    h_{2}(\mathbf{x}, \boldsymbol{\xi}) \\
    \vdots \\
    h_{r}(\mathbf{x}, \boldsymbol{\xi})
  \end{array}
  \right )\\
  \mbox{subject to}\\
  g_{j}(\mathbf{x}, \boldsymbol{\xi})\geq 0,  \ j=1, 2\dots, s, \\
  \end{array}
\end{equation}
where $\mathbf{x}$ is $k$-dimensional and $\boldsymbol{\xi}$ is $m$-dimensional. If
$\mathbf{x}$ or $\boldsymbol{\xi}$ are random,  then (\ref{eq0}) defines a multiobjective
stochastic optimisation problem,  see \citet{pr:95}. As shown in the following sections,  the
optimisation of a multiresponse surface can be proposed as a multiobjective stochastic
optimisation problem.

In this work,  the optimisation of a multiresponse surface is proposed as a multiobjective
stochastic optimisation problem. Section \ref{sec2} consider the notation and basic elements of
the multiresponse surface. Several previous approaches made are discussed in Section
\ref{sec3}. Section \ref{sec4} proposes the optimisation of a multiresponse surface as a
multiobjective stochastic optimisation problem and diverse solutions are proposed. Several
multiobjective stochastic solutions are studied in detail in Section \ref{sec5}. Finally,  a
real case from the literature is analysed in Section \ref{sec6}.

\section{Notation}\label{sec2}

A detailed discussion of multiresponse surface methodology may be found in \citet[Chap.
7]{kc:87} and \citet{kc:81}. For convenience,  the principal properties and usual notation is
restated here.

Let $N$ be the number of experimental runs and $r$ be the number of response variables which
can be measured for each setting of a group  of $n$ coded variables (also termed factors)
$x_{1},  x_{2}, \dots, x_{n}$. We assume that the response variables can be modelled by a
second order polynomial regression model in terms of $x_{j}^{Ôø?s}$. Hence,  the $k^{th}$
response model can be written as
\begin{equation}\label{lmu}
    \mathbf{Y}_{k} = \mathbf{X}_{k}\boldsymbol{\beta}_{k} + \boldsymbol{\varepsilon}_{k}
\end{equation}
where $\mathbf{Y}_{k}$ is an $N \times 1$ vector of observations on the $k^{th}$ response,
$\mathbf{X}_{k}$ is an $N \times p$ matrix of rank $p$ termed the design or regression matrix,
$p = 1 + n + n(n+1)/2$,  $\boldsymbol{\beta}_{k}$ is a $p \times 1$ vector of unknown constant
parameters,  and $\boldsymbol{\varepsilon}_{k}$ is a random error vector associated with the
$k^{th}$ response. In the present case, it is assumed that $\mathbf{X}_{1} = \cdots =
\mathbf{X}_{r} = \mathbf{X}$. Hence,  (\ref{lmu}) can be written as
\begin{equation}\label{lmm}
    \mathbf{Y} = \mathbf{X}\mathbb{B} + \mathbb{E}
\end{equation}
where $\mathbf{Y} = \left[\mathbf{Y}_{1}\vdots \mathbf{Y}_{2} \vdots \cdots \vdots
\mathbf{Y}_{r}\right]$,  $\mathbb{B} = \left[\boldsymbol{\beta}_{1}\vdots
\boldsymbol{\beta}_{2} \vdots \cdots \vdots \boldsymbol{\beta}_{r}\right]$ and $\mathbb{E} =
\left[\boldsymbol{\varepsilon}_{1}\vdots \boldsymbol{\varepsilon}_{2} \vdots \cdots \vdots
\boldsymbol{\varepsilon}_{r}\right]$,  such that $\mathbb{E} \sim \mathcal{N}_{N \times r}
(\mathbf{0},  \mathbf{I}_{N} \otimes \mathbf{\Sigma})$ i.e. $\mathbb{E}$ has an $N \times r$
matrix multivariate normal distribution with $\E (\mathbb{E}) = \mathbf{0}$ and $\Cov(\Vec
\mathbb{E}') = \mathbf{I}_{N} \otimes \mathbf{\Sigma}$,  where $\mathbf{\Sigma}$ is a $r \times
r$ positive definite matrix,  where if $\mathbf{A} = \left[\mathbf{A}_{1}\vdots
\mathbf{A}_{2}\vdots \cdots \vdots \mathbf{A}_{r}\right]$,  then $\Vec \mathbf{A} =
(\mathbf{A}'_{1},  \mathbf{A}'_{2},  \dots,  \mathbf{A}'_{r})'$ and $\otimes$ denotes the
direct (or Kronecker) product of matrices,  see \citet[Theorem 3.2.2,  p. 79]{mh:82}. In
addition let
\begin{description}
\item[$\centerdot$] $\mathbf{x} = \left(x_{1},  x_{2}, \dots, x_{n}\right)'$: The vector of controllable
variables or factors. Formally,  an $x_{i}$ variable is associated with each factor $A,  B,
...$

\item[$\centerdot$] $\widehat{\mathbb{B}} = \left[\widehat{\boldsymbol{\beta}}_{1}\vdots
\widehat{\boldsymbol{\beta}}_{2} \vdots \cdots \vdots \widehat{\boldsymbol{\beta}}_{r}\right]$:
The least squares estimator of $\mathbb{B}$ given by $\widehat{\mathbb{B}} =
(\mathbb{X}'\mathbb{X})^{-1}\mathbb{X}'\mathbf{Y}$,  from where
$\widehat{\boldsymbol{\beta}}_{k} = (\mathbb{X}' \mathbb{X})^{-1} \mathbb{X}' \mathbf{Y}_{k}$,
$k = 1,  2,  \dots,  r$. Moreover,  under the assumption that $\mathbb{E} \sim \mathcal{N}_{N
\times r} (\mathbf{0},  \mathbf{I}_{N} \otimes \mathbf{\Sigma})$,  then $\widehat{\mathbb{B}}
\sim \mathcal{N}_{p \times r} (\mathbb{B},  (\mathbf{X}'\mathbf{X})^{-1} \otimes
\mathbf{\Sigma})$,  with $\Cov(\Vec \widehat{\mathbb{B}}') = (\mathbf{X}'\mathbf{X})^{-1}
\otimes \mathbf{\Sigma}$.

\item[$\centerdot$] $\mathbf{z}(\mathbf{x}) = (1, x _{1},  x_{2},  \dots,  x_{n},  x _{1}^{2},  x_{2}^{2},
\dots,  x_{n}^{2},  x _{1} x_{2},  x_{1}x_{3}\dots,  x_{n-1}x_{n})'.$

\item[$\centerdot$] \hspace{-3cm} \vspace{-.85cm}
\begin{eqnarray*}
  \widehat{Y}_{k} (\mathbf{x}) &=& \mathbf{z}'(\mathbf{x})\widehat{\boldsymbol{\beta}}_{k} \\
&=& \widehat{\beta}_{0k} + \displaystyle\sum_{i=1}^{n}\widehat{\beta}_{ik}x_{i} +
\sum_{i=1}^{n}
    \widehat{\beta}_{iik}x_{i}^{2} +
    \sum_{i=1}^{n}\sum_{j>i}^{n}\widehat{\beta}_{ijk}x_{i}x_{j}:
\end{eqnarray*}
The response surface or predictor equation at the point $\mathbf{x}$ for the k$^{th}$ response
variable.

\item[$\centerdot$] $\widehat{\mathbf{Y}} (\mathbf{x}) = \left(\widehat{Y}_{1} (\mathbf{x}),  \widehat{Y}_{2}
(\mathbf{x}),  \dots,  \widehat{Y}_{r} (\mathbf{x})\right)' = \widehat{\mathbb{B}}'
\mathbf{z}(\mathbf{x})$: The multiresponse surface or predicted response vector at the point
$\mathbf{x}$.

\item[$\centerdot$]$\widehat{\mathbf{\Sigma}} = \displaystyle \frac{\mathbf{Y}'(\mathbf{I}_{N} -
\mathbf{X}(\mathbf{X}'\mathbf{X})^{-1}\mathbf{X}')\mathbf{Y}}{N-p}$: The estimator of the
variance-covariance matrix $\mathbf{\Sigma}$ such that $(N-p)\widehat{\mathbf{\Sigma}}$ has a
Wishart distribution with $(N-p)$ degrees of freedom and the parameter $\mathbf{\Sigma}$; this
fact is denoted as $(N-p)\widehat{\mathbf{\Sigma}} \sim \mathcal{W}_{r}(N-p, \mathbf{\Sigma})$.
Here,  $\mathbf{I}_{m}$ denotes an identity matrix of order $m$.

Finally,  note that
\begin{equation}\label{eq5}
    E(\widehat{\mathbf{Y}}(\mathbf{x})) =
    E(\widehat{\mathbb{B}}'\mathbf{z}(\mathbf{x})) = \mathbb{B}'\mathbf{z}(\mathbf{x})
\end{equation}
and
\begin{equation}\label{eq6}
    \Cov (\widehat{\mathbf{Y}}(\mathbf{x})) = \mathbf{z}'(\mathbf{x})(\mathbf{X}'
    \mathbf{X})^{-1}\mathbf{z}(\mathbf{x}) \mathbf{\Sigma}.
\end{equation}
An unbiased estimator of $\Cov (\widehat{\mathbf{Y}}(\mathbf{x}))$ is given by
\begin{equation}\label{eq7}
   \widehat{ \Cov} (\widehat{\mathbf{Y}}(\mathbf{x})) = \mathbf{z}'(\mathbf{x})(\mathbf{X}'
    \mathbf{X})^{-1}\mathbf{z}(\mathbf{x}) \widehat{\mathbf{\Sigma}}.
\end{equation}
\end{description}

\section{Multiresponse optimisation}\label{sec3}

In the following sections,  we make use of multiresponse optimisation and multiobjective (or
more general multicriteria) optimisation. For convenience,  the concepts and notation required
are listed below in terms of the estimated model of multiresponse surface optimisation.
Definitions and detailed properties may be found in \citet{kc:81},  \citet{kc:87},
\citet{rio:89},  \citet{s86},  \citet{m99},  \citet{v:72} and \citet{pr:95}.

The multiresponse optimisation (MRO) problem in general is proposed as
\begin{equation}\label{equiv1}
  \begin{array}{c}
    \build{\min}{}{\mathbf{x}}\widehat{\mathbf{Y}}(\mathbf{x}) = \
    \build{\min}{}{\mathbf{x}}
    \left(%
    \begin{array}{c}
      \widehat{Y}_{1} (\mathbf{x}) \\
      \widehat{Y}_{2} (\mathbf{x}) \\
      \vdots\\
      \widehat{Y}_{r} (\mathbf{x}) \\
    \end{array}%
    \right)\\
    \mbox{subject to}\\
    \mathbf{x} \in \mathfrak{X},
  \end{array}
\end{equation}
which is a deterministic nonlinear multiobjective optimisation problem,  see \citet{s86},
\citet{rio:89} and \citet{m99}; and where $\mathfrak{X}$ denotes the experimental region,
which in general is defined as a hypercube
$$
  \mathfrak{X}= \{\mathbf{x}|l_{i} < x_{i} < u_{i},  \quad i = 1,  2,  \dots,  n\},
$$
where $\mathbf{l} = \left(l_{1},  l_{2},  \dots,  l_{n}\right)'$,  defines the vector of lower
bounds of factors and $\mathbf{u} = \left(u_{1},  u_{2},  \dots,  u_{n}\right)'$,  define the
vector of upper bounds of factors. Alternatively,  the experimental region is defined as a
hypersphere
$$
  \mathfrak{X}= \{\mathbf{x}|\mathbf{x}'\mathbf{x} \leq c^{2},   c \in \Re\},
$$
where,  in general $c$ is determined by the experimental design model used,  see \citet{kc:87}.
Alternatively (\ref{equiv1}) can be written as
$$
  \build{\min}{}{\mathbf{x} \in \mathfrak{X}}\widehat{\mathbf{Y}}(\mathbf{x})
$$

In the response surface methodology context,  observe that in multiobjective optimisation
problems,  there rarely exists a point $\mathbf{x^{*}}$ which is considered as an optimum, i.e.
few cases satisfy the requirement that $ \widehat{Y}_{k} (\mathbf{x}) $ is minimum for all $k =
1, 2,  \dots,  r$. From the viewpoint of multiobjective optimisation,  this justifies the
following notion of the Pareto point,  which is more weakly defined than is an optimum point:
\begin{center}
\begin{minipage}[t]{4.5in}
\begin{it}
We say that $\widehat{\mathbf{Y}}^{*} (\mathbf{x})$ is a  \textit{Pareto point} of
$\widehat{\mathbf{Y}} (\mathbf{x})$,  if there is no other point $\widehat{\mathbf{Y}}^{1}
(\mathbf{x})$ such that $\widehat{\mathbf{Y}}^{1} (\mathbf{x}) \leq \widehat{\mathbf{Y}}^{*}
(\mathbf{x})$,  i.e. for all $k$,  $\widehat{Y}_{k}^{1} (\mathbf{x}) \leq \widehat{Y}_{k}^{*}
(\mathbf{x})$ and $\widehat{\mathbf{Y}}^{1} (\mathbf{x}) \neq \widehat{\mathbf{Y}}^{*}
(\mathbf{x})$.
\end{it}
\end{minipage}
\end{center}
Existence criteria for Pareto points in a multiobjective optimisation problem  and the
extension of scalar optimisation (\textit{Kuhn-Tucker's conditions}) to the vectorial case are
established in  \citet{s86},  \citet{rio:89} and \citet{m99}.

Methods for solving a multiobjective optimisation problem are based on the information
possessed about a particular problem. There are three possible scenarios: when the investigator
possesses either complete,  partial or null information,  see \citet{rio:89},  \citet{m99} and
\citet{s86}. In a response surface methodology context,  complete information means that the
investigator understands the population in such a way that it is possible to propose a
\textit{value function} reflecting the importance of each response variable,  where
\begin{center}
\begin{minipage}[t]{4.5in}
\begin{it}
    The value function is a function $f: \Re^{n} \rightarrow \Re$ such
    that $\min\widehat{\mathbf{Y}} (\mathbf{x^{*}}) < \min
    \widehat{\mathbf{Y}} (\mathbf{x}_{1}) \Leftrightarrow
    f(\widehat{\mathbf{Y}} (\mathbf{x^{*}})) <
    f(\widehat{\mathbf{Y}} (\mathbf{x}_{1})),  \quad \mathbf{x}^{*}\neq
    \mathbf{x}_{1}$.
\end{it}
\end{minipage}
\end{center}
In partial information,  the investigator knows the main response variable of the study very
well and this is sufficient support for the research. Finally,  under null information,  the
researcher only possesses information about the estimators of the response surface parameter,
and with this material an appropriate solution can be found too.

As can be observed,  all the approaches proposed in the literature are particular cases of the
models studied in multiobjective optimisation,  and in particular of the
\textit{$\epsilon$-constraint model} and the \textit{value function model} or a combination of
the two.Accordingly,  the equivalent nonlinear scalar optimisation problem of (\ref{equiv1}) is
of the form
\begin{equation}\label{eq1}
  \begin{array}{c}
    \build{\min}{}{\mathbf{x}}f\left(\widehat{\mathbf{Y}}(\mathbf{x})\right)\\
    \mbox{subject to} \\
    \mathbf{x} \in \mathfrak{X} \cap \mathfrak{S},
  \end{array}
\end{equation}
where $\mathfrak{S}$ is a subset generated by additional potential constraints,  which derive
from the particular technique used for establishing the equivalent deterministic scalar
optimisation problem (\ref{eq1}). In some particular cases of (\ref{eq1}),  a new fixed
parameter may appear,  such as a $\mathbf{w} = \left(w_{1},  w_{2}, \dots, w_{r}\right)'$,
vector of response weights and/or $\boldsymbol{\tau} = \left(\tau_{1}, \tau_{2}, \dots,
\tau_{r}\right)'$,  vector of target values for the response vector. Particular examples of
this equivalent univariate objective optimisation are the use of goal programming,  see
\citet{kea:08}, and of the $\epsilon$-constraint model,  see \citet{b:75},  among many others.

When uncertainty is assumed in an MRO problem,  in other words,  when the MRO problem is
considered as a stochastic program,  different approaches have been proposed,  see
\citet{n:08},  \citet{chh:01},  \citet{aea:08},  \citet{kc:81} and \citet{kc:87}. In particular
some of these approaches can be established as
\begin{equation}\label{eq8}
  \begin{array}{c}
    \build{\min}{}{\mathbf{x}}f\left(\widehat{\mathbf{Y}}(\mathbf{x})\right)\\
    \mbox{subject to} \\
    \mathbf{x} \in \mathfrak{X} \cap \mathfrak{S} \\
    \widehat{\mathbb{B}} \sim \mathcal{N}_{p \times r} (\mathbb{B},
    (\mathbf{X}'\mathbf{X})^{-1} \otimes \mathbf{\Sigma}) \\
    (N-p)\widehat{\mathbf{\Sigma}}\sim \mathcal{W}_{r}(N-p,  \mathbf{\Sigma}),
  \end{array}
\end{equation}
where $\widehat{\mathbb{B}}$ and $\widehat{\mathbf{\Sigma}}$ are independent.

In addition,  it is sometimes assumed that $\mathbf{w}$ and/or $\boldsymbol{\tau}$ are
stochastic,  and diverse strategies have been proposed to obtain particular values for these,
including Group Decision Making,  among others,  see \citet{kc:81} and \citet{hbdn:10}.

In general terms,  thus,  the MRO problem under uncertainty has been addressed as follows:
\begin{enumerate}
  \item It is considered as a deterministic multiobjective optimisation problem (\ref{equiv1}).
  \item An equivalent deterministic univariate optimisation problem, such as goal programming,  is proposed (\ref{eq1}).
  \item In the equivalent deterministic univariate optimisation problem,
uncertainty is assumed (\ref{eq8}).
\end{enumerate}

\section{Proposed approach}\label{sec4}

In the univariate case,  \citet{dea:05} considered the problem as a stochastic optimisation
program. The approach proposed in the present paper consists in extending this idea to
multiresponse optimisation. Specifically,  we propose the MRO as the following nonlinear
multiobjective stochastic optimisation problem from the beginning:
\begin{equation}\label{equiv2}
  \begin{array}{c}
    \build{\min}{}{\mathbf{x}}\widehat{\mathbf{Y}}\left(\mathbf{x}, \widehat{\mathbb{B}}\right)\\
    \mbox{subject to}\\
    \mathbf{x} \in \mathfrak{X}\\
    \widehat{\mathbb{B}} \sim \mathcal{N}_{p \times r} (\mathbb{B},
    (\mathbf{X}'\mathbf{X})^{-1} \otimes \mathbf{\Sigma}) \\
    (N-p)\widehat{\mathbf{\Sigma}}\sim \mathcal{W}_{r}(N-p,  \mathbf{\Sigma}).
  \end{array}
\end{equation}
where $ \widehat{\mathbf{Y}}(\mathbf{x},  \widehat{\mathbb{B}}) \equiv
\widehat{\mathbf{Y}}(\mathbf{x})$ and $\widehat{Y}_{k}(\mathbf{x},
\widehat{\boldsymbol{\beta}}_{k}) \equiv \widehat{Y}_{k}(\mathbf{x})$.

The solution of (\ref{equiv2}) can be applied to any model (technique,  method or solution)
under multiobjective stochastic optimisation,  which in general,  is a multidimensional
extension of stochastic optimisation models,  see \citet{v:72},  \citet{dea:05} and
\citet{pr:95}.

\subsection{Multiobjective stochastic optimisation approaches}\label{sub41}

In this subsection we propose (\ref{equiv2}) under diverse multiobjective stochastic
optimisation approaches. The properties of the solution obtained under the different approaches
are described in detail by \citet{k:63},  \citet{sm:84} and \citet{pr:95}.

As shown below,  each multiobjective stochastic optimisation approach can be proposed in
several ways. In some cases,  this possibility is a consequence of assuming that the response
variables are correlated or not.

\subsubsection{Multiobjective expected value solution,  multiobjective E-model}

Point $\mathbf{x} \in \mathfrak{X}$ is the expected value solution to (\ref{equiv2}) if it is
an efficient solution in the Pareto sense to the following deterministic multiobjective
optimisation problem
\begin{equation}\label{solE}
    \build{\min}{}{\mathbf{x} \in \mathfrak{X}}\E\left(\widehat{\mathbf{Y}}\left(\mathbf{x},
    \widehat{\mathbb{B}}\right)\right)
\end{equation}

\subsubsection{Multiobjective minimum variance solution,  multiobjective V-model}

The $\mathbf{x} \in \mathfrak{X}$ point is the minimum variance solution to problem
(\ref{equiv2}) if it is an efficient solution in the Pareto sense of the deterministic
multiobjective optimisation problem
\begin{equation}\label{solV1}
    \build{\min}{}{\mathbf{x} \in \mathfrak{X}}
    \left(%
    \begin{array}{c}
      \Var\left(\widehat{Y}_{1} (\mathbf{x})\right) \\
      \Var\left(\widehat{Y}_{2} (\mathbf{x})\right) \\
      \vdots\\
      \Var\left(\widehat{Y}_{r} (\mathbf{x})\right). \\
    \end{array}%
    \right)
\end{equation}
This efficient solution is adequate if it is assumed that the response variables are
Uncorrelated. However,  if the response variables are assumed to be correlated a better one is:
\begin{equation}\label{solV2}
    \build{\min}{}{\mathbf{x} \in \mathfrak{X}}\Cov\left(\widehat{\mathbf{Y}}\left(\mathbf{x},
    \widehat{\mathbb{B}}\right)\right)
\end{equation}

\subsubsection{Multiobjective expected value standard deviation solution,  multiobjective modified E-model}

Point $\mathbf{x} \in \mathfrak{X}$ is an expected value standard deviation solution to the
problem (\ref{equiv2}) if it is an efficient solution in the Pareto sense of the mixed
deterministic multiobjective-matrix optimisation problem
\begin{equation}\label{solV}
    \build{\min}{}{\mathbf{x} \in \mathfrak{X}}
    \left[
    \begin{array}{c}
      \E\left(\widehat{\mathbf{Y}}\left(\mathbf{x}, \widehat{\mathbb{B}}\right)\right) \\
      \left(\Cov\left(\widehat{\mathbf{Y}}\left(\mathbf{x}, \widehat{\mathbb{B}}\right)\right)\right)^{1/2}
    \end{array}
    \right ]
\end{equation}
where $\left(\mathbf{A}^{1/2}\right)^{2} = \mathbf{A}$,  see \citet[Appendix]{mh:82}.

We now define the concept of the efficient solution of multiobjective minimum risk of joint
aspiration level $\boldsymbol{\tau} = (\tau_{1},  \tau_{2},  \dots,  \tau_{r})'$ and the
efficient solution with a joint probability $\alpha$. The two solutions are obtained by
applying the multivariate versions of minimum risk and the Kataoka criteria,  respectively,
referred to in the literature as criteria of maximum probability or satisfying criteria,  due
to the fact that,  as shown below,  in both cases the criteria to be used provide,  in one way
or another,  ``good" solutions in terms of probability,  see \citet{k:63}.

\subsubsection{Multiobjective minimum risk solution of joint aspiration level $\boldsymbol{\tau}$,
multiobjective modified $\P$-model}

Point $\mathbf{x} \in \mathfrak{X}$ is a minimum risk solution of joint aspiration level
$\boldsymbol{\tau}$ to problem (\ref{equiv2}) if it constitutes an efficient solution in the
Pareto sense of the multiobjective stochastic optimisation problem
\begin{equation}\label{solMR1}
    \build{\min}{}{\mathbf{x} \in \mathfrak{X}}
    \left(%
    \begin{array}{c}
      \P\left(\widehat{Y}_{1} (\mathbf{x}) \leq \tau_{1}\right) \\
      \P\left(\widehat{Y}_{2} (\mathbf{x}) \leq \tau_{2}\right) \\
      \vdots\\
      \P\left(\widehat{Y}_{r} (\mathbf{x}) \leq \tau_{r}\right) \\
    \end{array}%
    \right).
\end{equation}
It is also possible to consider the following alternative multiobjective $\P$-model
\begin{equation}\label{solMR2}
    \build{\min}{}{\mathbf{x} \in \mathfrak{X}}
     \P\left(%
    \begin{array}{c}
      \widehat{Y}_{1} (\mathbf{x}) \leq \tau_{1} \\
      \widehat{Y}_{2} (\mathbf{x}) \leq \tau_{2} \\
      \vdots\\
      \widehat{Y}_{r} (\mathbf{x}) \leq \tau_{r} \\
    \end{array}%
    \right).
\end{equation}
Again,  (\ref{solMR2}) is more adequate if the response variables are correlated. However
(\ref{solMR2}) is considerably more complicated to solve than (\ref{solMR1}). When $r = 2$,
\citet{pr:70} proposed an algorithm for a similar problem (probabilistic constrained
programming),  which can be applied to solve (\ref{solMR2}).

\subsubsection{Multiobjective Kataoka solution with probability $\alpha$}

Point $\mathbf{x} \in \mathfrak{X}$ is a multiobjective Kataoka solution with probability
$\alpha$ (fixed) to problem (\ref{equiv2}) if it is an efficient solution in the Pareto sense
of the multiobjective optimisation problem
\begin{equation}\label{solK1}
  \begin{array}{c}
    \build{\min}{}{\mathbf{x},  \boldsymbol{\tau}} \boldsymbol{\tau}\\
    \mbox{subject to}\\
      \P\left(\widehat{Y}_{k} (\mathbf{x}) \leq \tau_{k}\right) = \alpha,  \ k = 1, 2,  \dots,  r \\
    \mathbf{x} \in \mathfrak{X}.
  \end{array}
\end{equation}
Alternatively (\ref{solK1}) can be proposed as
\begin{equation}\label{solK2}
  \begin{array}{c}
    \build{\min}{}{\mathbf{x},  \boldsymbol{\tau}} \boldsymbol{\tau}\\
    \mbox{subject to}\\
      \P\left(%
    \begin{array}{c}
      \widehat{Y}_{1} (\mathbf{x}) \leq \tau_{1} \\
      \widehat{Y}_{2} (\mathbf{x}) \leq \tau_{2} \\
      \vdots\\
      \widehat{Y}_{r} (\mathbf{x}) \leq \tau_{r} \\
    \end{array}%
    \right) = \alpha\\
    \mathbf{x} \in \mathfrak{X},
  \end{array}
\end{equation}
Note that (\ref{solK1}) and (\ref{solK2}) are multiobjective probabilistic constrained
programming,  see \citet{chc:63},  \citet{sm:84} and \citet{pr:95}.

Many other approaches can be used to solve (\ref{equiv2}). For example,  \citet{sm:84} proposed
a stochastic version of the sequential technique,  termed the Lexicographic method,  for
solving (\ref{solV1}) and (\ref{solMR1}) or direct application to (\ref{equiv2}); among many
other options.

In all cases,  observe that
$$
  \E\left(\widehat{\mathbf{Y}}\left(\mathbf{x}, \widehat{\mathbb{B}}\right)\right) =
  \mathbf{Y}\left(\mathbf{x}, \mathbb{B}\right) \mbox{ and }
  \Cov\left(\widehat{\mathbf{Y}}\left(\mathbf{x}, \widehat{\mathbb{B}}\right)\right) =
  \mathbf{z}'(\mathbf{x})(\mathbf{X}' \mathbf{X})^{-1}\mathbf{z}(\mathbf{x}) \mathbf{\Sigma},
$$
in general are unknown. Then,  from a practical point of view,  and having the final expression
of the equivalent deterministic problem of (\ref{equiv2}),  $
\E\left(\widehat{\mathbf{Y}}\left(\mathbf{x}, \widehat{\mathbb{B}}\right)\right)$ and
$\Cov\left(\widehat{\mathbf{Y}}\left(\mathbf{x}, \widehat{\mathbb{B}}\right)\right)$ should be
replaced by their corresponding estimators
$$
  \E\left(\widehat{\mathbf{Y}}\left(\mathbf{x}, \widehat{\mathbb{B}}\right)\right) =
  \widehat{\mathbf{Y}}\left(\mathbf{x}, \widehat{\mathbb{B}}\right) \mbox{ and }
  \widehat{\Cov}\left(\widehat{\mathbf{Y}}\left(\mathbf{x}, \widehat{\mathbb{B}}\right)\right) =
  \mathbf{z}'(\mathbf{x})(\mathbf{X}' \mathbf{X})^{-1}\mathbf{z}(\mathbf{x})
  \widehat{\mathbf{\Sigma}}.
$$

\section{Equivalent deterministic programs }\label{sec5}

In this section we study  several particular equivalent deterministic programs from
(\ref{solV2}) in detail.

\subsection{Multiobjective $V$-model}

Taking into account the final comment in Section \ref{sec4},  our intention is to solve the
matrix optimisation problem
\begin{equation}\label{solV21}
    \build{\min}{}{\mathbf{x} \in \mathfrak{X}}\widehat{\Cov}\left(\widehat{\mathbf{Y}}\left(\mathbf{x},
    \widehat{\mathbb{B}}\right)\right).
\end{equation}
For the sake of convenience,  in this section we denote
$\widehat{\Cov}\left(\widehat{\mathbf{Y}}\left(\mathbf{x},  \widehat{\mathbb{B}}\right)\right)$
as $\widehat{\Cov}\left(\widehat{\mathbf{Y}}\left(\mathbf{x} \right)\right)$. Obviously,  the
difficulty in expressing the problem in this way lies in defining the meaning of the minimum of
a matrix function. The idea of minimising a matrix function,  and in particular a matrix of
variance-covariance,  has been studied with respect to various areas of statistical theory. For
example,  when regression estimators are determined for a multivariate general linear model,
this is done by minimising the determinant or the trace of sums of squares and cross-products
matrix of the error,  see \citet{g:77}. Similarly,  the choice or comparison of experimental
design models is done by minimising a function of the variance-covariance matrix of treatment
estimators,  see \citet{kc:87} and \citet{ad:97}.

Fortunately,  it is possible to reduce the nonlinear matrix minimisation problem (\ref{solV21})
to a univariate nonlinear minimisation problem by taking into account the following
considerations. Observe that the procedure described here is just one of various possible
options,  see \citet{rio:89} and \citet{m99}.

Assume that $\widehat{\Cov}\left(\widehat{\mathbf{Y}}\left(\mathbf{x} \right)\right)$ is a
positive definite matrix for all $\mathbf{x}$,  denoting it as
$\widehat{\Cov}\left(\widehat{\mathbf{Y}}\left(\mathbf{x} \right)\right) > \mathbf{0}$. Now,
let $\mathbf{x}_{1}$ and $\mathbf{x}_{2}$ be two possible values of the vector $\mathbf{x}$ and
let $\mathbf{B} = \widehat{\Cov}\left(\widehat{\mathbf{Y}}\left(\mathbf{x}_{1} \right)\right) -
\widehat{\Cov}\left(\widehat{\mathbf{Y}}\left(\mathbf{x}_{2} \right)\right)$. Then we say that
\begin{equation}\label{nd}
    \widehat{\Cov}\left(\widehat{\mathbf{Y}}\left(\mathbf{x}_{1} \right)\right) <
    \widehat{\Cov}\left(\widehat{\mathbf{Y}}\left(\mathbf{x}_{2} \right)\right) \ \Leftrightarrow \
    \mathbf{B} < \mathbf{0},
\end{equation}
i.e. if the matrix $\mathbf{B}$ is a negative definite matrix. Moreover,  note that
$\widehat{\Cov}\left(\widehat{\mathbf{Y}}\left(\mathbf{x}_{1} \right)\right)$ and $
\widehat{\Cov}\left(\widehat{\mathbf{Y}}\left(\mathbf{x}_{2} \right)\right)$,  are
diagonalizable. Then,  let $D_{\mathbf{x}_{1}}$ and $D_{\mathbf{x}_{2}}$ be the diagonal
matrixes associated with $\widehat{\Cov}\left(\widehat{\mathbf{Y}}\left(\mathbf{x}_{1}
\right)\right)$ and $ \widehat{\Cov}\left(\widehat{\mathbf{Y}}\left(\mathbf{x}_{2}
\right)\right)$,  respectively; with $D_{\mathbf{x}_{1}} = \diag(\alpha_{1}, \dots,
\alpha_{r})$,  $\alpha_{1}> \cdots
> \alpha_{r} > 0$ and $D_{\mathbf{x}_{2}} = \diag(\gamma_{1},  \dots, \gamma_{r})$,
$\gamma_{1}> \cdots > \gamma_{r} > 0$,  where $\alpha_{j}$ and $\gamma_{j}$ denote the
eigenvalues of $\widehat{\Cov}\left(\widehat{\mathbf{Y}}\left(\mathbf{x}_{1} \right)\right)$
and $ \widehat{\Cov}\left(\widehat{\mathbf{Y}}\left(\mathbf{x}_{2} \right)\right)$,
respectively. Thus,  expression (\ref{nd}) can alternatively be presented as:
$$
  \widehat{\Cov}\left(\widehat{\mathbf{Y}}\left(\mathbf{x}_{1} \right)\right) <
  \widehat{\Cov}\left(\widehat{\mathbf{Y}}\left(\mathbf{x}_{2} \right)\right) \ \Leftrightarrow \
  D_{\mathbf{x}_{1}} -D_{\mathbf{x}_{2}}< \mathbf{0},
$$
i.e.
\begin{equation}\label{pareto}
    \widehat{\Cov}\left(\widehat{\mathbf{Y}}\left(\mathbf{x}_{1} \right)\right) <
    \widehat{\Cov}\left(\widehat{\mathbf{Y}}\left(\mathbf{x}_{2} \right)\right) \ \Leftrightarrow \
    \build{\alpha_{j} - \gamma_{j }}{}{j = 1,  \dots,  r} \ < 0 ,
\end{equation}
and $\widehat{\Cov}\left(\widehat{\mathbf{Y}}\left(\mathbf{x}_{1} \right)\right) \neq
\widehat{\Cov}\left(\widehat{\mathbf{Y}}\left(\mathbf{x}_{2} \right)\right)$; which defines a
weak Pareto order,  see \citet{s86},  \citet{rio:89} and \citet{m99}. Then from \citet{s86},
\citet{rio:89} and \citet{m99},  there exists a function $g: \mathcal{S}\rightarrow \Re$,  such
that
\begin{equation}\label{citerio}
    \widehat{\Cov}\left(\widehat{\mathbf{Y}}\left(\mathbf{x}_{1} \right)\right) <
    \widehat{\Cov}\left(\widehat{\mathbf{Y}}\left(\mathbf{x}_{2} \right)\right)
    \Leftrightarrow
    g\left(\widehat{\Cov}\left(\widehat{\mathbf{Y}}\left(\mathbf{x}_{1} \right)\right)\right)
< g\left(\widehat{\Cov}\left(\widehat{\mathbf{Y}}\left(\mathbf{x}_{2} \right)\right)\right).
\end{equation}
where $\widehat{\Cov}\left(\widehat{\mathbf{Y}}\left(\mathbf{x} \right)\right) \in
\mathcal{S}\subset \Re^{r(r+1)/2}$ and $\mathcal{S}$ is the set of positive definite matrices.
From (\ref{citerio}),  \citet{s86},  \citet{rio:89} and \citet{m99} prove that the non-linear
matrix minimisation problem (\ref{solV21}) is reduced in the following scalar non-linear
minimisation problem
\begin{equation}\label{equi1}
    \build{\min}{}{\mathbf{x} \in \mathfrak{X}}g\left(\widehat{\Cov}\left(\widehat{\mathbf{Y}}\left(\mathbf{x},
    \widehat{\mathbb{B}}\right)\right)\right).
\end{equation}

Unfortunately or otherwise,  the function $g(\cdot)$ is not unique. For example,  in other
statistical contexts we can find the following commonly used functions $g(\cdot)$,  see
\citet{g:77}:
\begin{enumerate}
    \item The trace of the matrix
    $\widehat{\Cov}\left(\widehat{\mathbf{Y}}\left(\mathbf{x} \right)\right)$;
    \begin{eqnarray*}
      g\left(\widehat{\Cov}\left(\widehat{\mathbf{Y}}\left(\mathbf{x} \right)\right)\right) &=& \tr
        \left(\widehat{\Cov}\left(\widehat{\mathbf{Y}}\left(\mathbf{x} \right)\right)\right) \\
&=& \mathbf{z}'(\mathbf{x})(\mathbf{X}' \mathbf{X})^{-1}\mathbf{z}(\mathbf{x})
        \displaystyle\sum_{j}^{r} \widehat{\sigma}_{jj}.
    \end{eqnarray*}
    \item The determinant of the matrix
    $\widehat{\Cov}\left(\widehat{\mathbf{Y}}\left(\mathbf{x} \right)\right)$;
    \begin{eqnarray*}
      g\left(\widehat{\Cov}\left(\widehat{\mathbf{Y}}\left(\mathbf{x} \right)\right)\right) &=&
        \left|\widehat{\Cov}\left(\widehat{\mathbf{Y}}\left(\mathbf{x} \right)\right)\right| \\
&=& \left[\mathbf{z}'(\mathbf{x})(\mathbf{X}' \mathbf{X})^{-1}\mathbf{z}(\mathbf{x})\right]^{r}
        |\widehat{\mathbf{\Sigma}}|
    \end{eqnarray*}
    \item The sum of all the elements of the matrix  $\widehat{\Cov}\left(\widehat{\mathbf{Y}}(\mathbf{x})\right)$;
        $$g\left(\widehat{\Cov}\left(\widehat{\mathbf{Y}}\left(\mathbf{x} \right)\right)\right) =
            \mathbf{z}'(\mathbf{x})(\mathbf{X}' \mathbf{X})^{-1}\mathbf{z}(\mathbf{x})
            \displaystyle\sum_{j, k =1}^{r} \widehat{\sigma}_{jk}.
        $$
    \item $g\left(\widehat{\Cov}\left(\widehat{\mathbf{Y}}\left(\mathbf{x} \right)\right)\right)=
        \lambda_{\max}\left(\widehat{\Cov}\left(\widehat{\mathbf{Y}}\left(\mathbf{x} \right)\right)\right)$,
        where $\lambda_{\max}$ is the maximum eigenvalue
        of the matrix $\widehat{\Cov}\left(\widehat{\mathbf{Y}}\left(\mathbf{x} \right)\right)$.
    \item $g\left(\widehat{\Cov}\left(\widehat{\mathbf{Y}}\left(\mathbf{x} \right)\right)\right)=
        \lambda_{\min}\left(\widehat{\Cov}\left(\widehat{\mathbf{Y}}\left(\mathbf{x} \right)\right)\right)$,
        where $\lambda_{\min}$ is the minimum eigenvalue of
        the matrix  $\widehat{\Cov}\left(\widehat{\mathbf{Y}}\left(\mathbf{x} \right)\right)$.
    \item $g\left(\widehat{\Cov}\left(\widehat{\mathbf{Y}}\left(\mathbf{x} \right)\right)\right)=
        \lambda_{j}\left(\widehat{\Cov}\left(\widehat{\mathbf{Y}}\left(\mathbf{x} \right)\right)\right)$,
        where $\lambda_{j}$ is the $j$-th eigenvalue of
        the matrix $\widehat{\Cov}\left(\widehat{\mathbf{Y}}\left(\mathbf{x} \right)\right)$,  among others.
\end{enumerate}

But observe that
\begin{eqnarray*}
  \lambda_{j}\left(\widehat{\Cov}\left(\widehat{\mathbf{Y}}\left(\mathbf{x}
    \right)\right)\right) &=& \lambda_{j}\left(\mathbf{z}'(\mathbf{x})(\mathbf{X}' \mathbf{X})^{-1}\mathbf{z}(\mathbf{x})
    \widehat{\mathbf{\Sigma}}\right) \\
&=& \mathbf{z}'(\mathbf{x})(\mathbf{X}' \mathbf{X})^{-1}\mathbf{z}(\mathbf{x})\lambda_{j}\left(
   \widehat{\mathbf{\Sigma}}\right).
\end{eqnarray*}
Hence we can conclude that as a consequence of the structure of the covariance matrix
$\widehat{\Cov}\left(\widehat{\mathbf{Y}}\left(\mathbf{x}\right)\right)$,  for all the
particular definitions of the function $g$ considered above,  the scalar non-linear
minimisation problem (\ref{equi1}) has a unique solution given by the solution of the
non-linear minimisation problem
\begin{equation}\label{equi2}
    \build{\min}{}{\mathbf{x} \in \mathfrak{X}}\mathbf{z}'(\mathbf{x})(\mathbf{X}' \mathbf{X})^{-1}\mathbf{z}(\mathbf{x}).
\end{equation}

\subsection{Multiobjective $\P$-model}

Proceeding as in \citet{dea:05},  the equivalent multiobjective deterministic problem to
(\ref{equiv2}) via the $\P$-model (\ref{solMR1}) is
\begin{equation}\label{dsolMR1}
    \build{\min}{}{\mathbf{x} \in \mathfrak{X}}
    \left(%
    \begin{array}{c}\displaystyle
      \frac{\tau_{1} - z'(\mathbf{x})\widehat{\boldsymbol{\beta}}_{1}}{
        \sqrt{\widehat{\sigma}_{11}\ z'(\mathbf{x})(\mathbf{X}'\mathbf{X})^{-1}z(\mathbf{x})}} \\
      \displaystyle
      \frac{\tau_{2} - z'(\mathbf{x})\widehat{\boldsymbol{\beta}}_{2}}{
        \sqrt{\widehat{\sigma}_{22}\ z'(\mathbf{x})(\mathbf{X}'\mathbf{X})^{-1}z(\mathbf{x})}} \\
      \vdots\\
      \displaystyle
      \frac{\tau_{r} - z'(\mathbf{x})\widehat{\boldsymbol{\beta}}_{r}}{
        \sqrt{\widehat{\sigma}_{rr}\ z'(\mathbf{x})(\mathbf{X}'\mathbf{X})^{-1}z(\mathbf{x})}} \\
    \end{array}%
    \right).
\end{equation}

\subsection{Multiobjective Kataoka model}

From \citet{dea:05},  the equivalent multiobjective deterministic problem to (\ref{equiv2}) via
the Kataoka model (\ref{solK1}) is given by
\begin{equation}\label{dsolMR2}
    \build{\min}{}{\mathbf{x} \in \mathfrak{X}}
    \left(%
    \begin{array}{c}
      z'(\mathbf{x}) \widehat{\boldsymbol{\beta}_{1}} + \Phi^{-1}(\delta) \
       \sqrt{\widehat{\sigma}_{11}\ z'(\mathbf{x})(\mathbf{X}'\mathbf{X})^{-1}z(\mathbf{x})}\\
      z'(\mathbf{x}) \widehat{\boldsymbol{\beta}_{2}} + \Phi^{-1}(\delta) \
       \sqrt{\widehat{\sigma}_{22}\ z'(\mathbf{x})(\mathbf{X}'\mathbf{X})^{-1}z(\mathbf{x})}\\
      \vdots\\
      z'(\mathbf{x}) \widehat{\boldsymbol{\beta}_{r}} + \Phi^{-1}(\delta) \
       \sqrt{\widehat{\sigma}_{rr}\ z'(\mathbf{x})(\mathbf{X}'\mathbf{X})^{-1}z(\mathbf{x})}\\
    \end{array}%
    \right).
\end{equation}
where $\Phi$ denotes the distribution function of the standard Normal distribution.

Similar equivalent multiobjective deterministic problems to (\ref{equiv2}) are obtained by
applying the other stochastic solutions described in Section \ref{sec4}. Note that,  if each
stochastic solution is combined with each multiobjective optimisation technique,  an infinite
number of possible solutions to (\ref{equiv2})is obtained. For example,  note that the function
of value $f(\cdot)$ may take an infinite number of forms. One of these particular forms is the
weighting method. Under this approach,  problem (\ref{dsolMR2}) can be restated as:
\begin{equation}\label{WsolMR2}
    \build{\min}{}{\mathbf{x} \in \mathfrak{X}}
      \sum_{k=1}^{r} w_{k}\left\{z'(\mathbf{x}) \widehat{\boldsymbol{\beta}_{k}} + \Phi^{-1}(\delta) \
      \sqrt{\widehat{\sigma}_{kk}\ z'(\mathbf{x})(\mathbf{X}'\mathbf{X})^{-1}z(\mathbf{x})}\right\}\\
\end{equation}
such that $\sum_{k=1}^{r}w_{k} = 1$,  $w_{k} \geq 0$ $\forall$ $k= 1,  2,  \dots,  r$: where
$w_{k}$ weights the importance of each characteristic. The solution $\mathbf{x} \in
\mathfrak{X}$ of (\ref{WsolMR2}) can be termed the multiobjective Kataoka solution with
probability $\alpha$ to problem (\ref{equiv2}),  via the weighting method.

Similarly,  $\mathbf{x} \in \mathfrak{X}$ is the multiobjective Kataoka solution with
probability $\alpha$ to the problem (\ref{equiv2}),  via goal programming if $\mathbf{x} \in
\mathfrak{X}$ is
\begin{equation}\label{GPsolMR2}
  \begin{array}{c}
    \build{\min}{}{\mathbf{x} \in \mathfrak{X}}\displaystyle\sum_{k=1}^{p}w_{k}(d_{k}^{+} + d_{k}^{-}) \\
    \mbox{subject to} \\
    z'(\mathbf{x}) \widehat{\boldsymbol{\beta}_{k}} + \Phi^{-1}(\delta) \
    \sqrt{\widehat{\sigma}_{kk}\ z'(\mathbf{x})(\mathbf{X}'\mathbf{X})^{-1}
    z(\mathbf{x})} - d_{k}^{+} + d_{k}^{-} = \tau_{k},  \ k = 1, 2, \dots, r,
  \end{array}
\end{equation}
where
\begin{eqnarray*}
  d_{k}^{+} &=& \frac{1}{2} \left( \left|z'(\mathbf{x}) \widehat{\boldsymbol{\beta}_{k}} + \Phi^{-1}(\delta) \
        \sqrt{\widehat{\sigma}_{kk}\ z'(\mathbf{x})(\mathbf{X}'\mathbf{X})^{-1}z(\mathbf{x})} - \tau_{k}
        \right| \right .\\
&& \qquad\left. + \left(z'(\mathbf{x}) \widehat{\boldsymbol{\beta}_{k}} + \Phi^{-1}(\delta) \
        \sqrt{\widehat{\sigma}_{kk}\ z'(\mathbf{x})(\mathbf{X}'\mathbf{X})^{-1}z(\mathbf{x})} -
        \tau_{k}\right)\right),  \\
  d_{k}^{-} &=& \frac{1}{2} \left( \left|z'(\mathbf{x}) \widehat{\boldsymbol{\beta}_{k}} + \Phi^{-1}(\delta) \
        \sqrt{\widehat{\sigma}_{kk}\ z'(\mathbf{x})(\mathbf{X}'\mathbf{X})^{-1}z(\mathbf{x})} -
        \tau_{k} \right| \right .\\
&& \qquad\left. - \left(z'(\mathbf{x}) \widehat{\boldsymbol{\beta}_{k}} + \Phi^{-1}(\delta) \
        \sqrt{\widehat{\sigma}_{kk}\ z'(\mathbf{x})(\mathbf{X}'\mathbf{X})^{-1}z(\mathbf{x})} -
        \tau_{k}\right)\right).
\end{eqnarray*}

\section{Application}\label{sec6}

A real case from the literature is analysed by applying several approaches and particular
solutions from multiobjective stochastic optimisation.

In the following example,  taken from \citet{pi:93},  there are two response variables
$\mathbf{Y} = (Y_{1},  Y_{2})'$ with 4 replicates and three setting variables $\mathbf{x} =
(x_{1},  x_{2},  x_{3})'$. The experimental data are shown in Table \ref{table3}. It is assumed
that the targets of the responses are $\boldsymbol{\tau} = (\tau_{1},  \tau_{2})'=(103,  73)'$.

\begin{table}[h]
  \centering
  \caption{Experimental data for the numerical example.}\label{table3}
  \begin{scriptsize}
  \begin{minipage}[t]{460pt}
  \begin{tabular}{cccc|cccc|cccc}
    \hline\hline
    \multicolumn{4}{c|}{Replicate} & 1 & 2 & 3 & 4 & 1 & 2 & 3 & 4\\
    \hline
    ID & $x_{1}$ & $x_{2}$ & $x_{3}$ & \multicolumn{4}{c|}{$Y_{1}$} &  \multicolumn{4}{|c}{$Y_{2}$} \\
    \hline\hline
    8 & 1 & 1 & 1 & 104.45 & 105.03 & 99.79 & 104.92 & 76.90 & 77.03 & 67.99 & 75.77 \\
    4 & 1 & 1 & -1 & 104.12 & 104.80 & 104.20 & 104.34 & 72.99 & 74.25 & 73.94 & 73.28 \\
    6 & 1 & -1 & 1 & 98.73 & 99.36 & 102.84 & 94.24 & 67.10 & 63.61 & 68.65 & 62.42 \\
    2 & 1 & -1 & -1 & 100.19 & 99.63 & 100.27 & 100.60 & 67.03 & 66.18 & 66.58 & 67.94 \\
    7 & -1 & 1 & 1 & 103.15 & 106.96 & 107.62 & 103.44 & 71.68 & 76.27 & 77.50 & 76.37 \\
    3 & -1 & 1 & -1 & 106.08 & 105.64 & 105.67 & 105.39 & 72.94 & 72.85 & 72.58 & 72.38 \\
    5 & -1 & -1 & 1 & 113.52 & 111.12 & 112.85 & 106.67 & 68.29 & 68.47 & 68.96 & 64.71 \\
    1 & -1 & -1 & -1 & 109.90 & 109.76 & 110.70 & 109.77 & 67.70 & 67.24 & 67.96 & 66.93 \\
    \hline\hline
  \end{tabular}
  \end{minipage}
  \end{scriptsize}
\end{table}

Equations (\ref{eq13}) and (\ref{eq14}) are response surfaces for $Y_{1}$ and $Y_{2}$.
\begin{eqnarray}\label{eq13}
  \widehat{Y}_{1}(\mathbf{x}) &=& 104.86 - 3.147x_{1} - 0.142x_{2} - 0.199x_{3} + 2.379x_{1}x_{2} \nonumber \\
&& \hspace{6cm} - 0.35x_{1}x_{3} - 0.106 x_{2}x_{3}\\
  \label{eq14}
  \widehat{Y}_{2}(\mathbf{x}) &=& 70.45 - 0.348x_{1} + 3.59x_{2} + 0.28x_{3} + 0.323x_{1}x_{2} \nonumber \\
&& \hspace{6cm}  - 0.45x_{1}x_{3} + 0.614 x_{2}x_{3}
\end{eqnarray}
From which the multiresponse optimisation problem is given as
\begin{equation}\label{eq16}
    \build{\min}{}{\mathbf{x}\in \mathfrak{X}}\widehat{\mathbf{Y}}(\mathbf{x}) = \
    \build{\min}{}{\mathbf{x}\in \mathfrak{X}}
    \left(%
    \begin{array}{c}
      \widehat{Y}_{1} (\mathbf{x}) \\
      \widehat{Y}_{2} (\mathbf{x}) \\
    \end{array}%
    \right)
\end{equation}
where $\mathfrak{X} = \{\mathbf{x}|x_{i} \in [-1, 1],  \ i = 1, 2, 3\}$.

Now,  assume that the importance of each response variable must be assessed from the decision
makers' viewpoint. For the purposes of this example,  consider $\mathbf{w}=(w_{1}, w_{2})' =
(0.285,  0.715)'$.

From (\ref{eq13}) and (\ref{eq14})%
{\footnotesize $$
  \hspace{3.2cm}
  \begin{array}{ccccccc}
    \quad\widehat{\boldsymbol{\beta}}_{0}\quad& \quad\widehat{\boldsymbol{\beta}}_{1}\quad
& \quad\widehat{\boldsymbol{\beta}}_{2}\quad & \quad\widehat{\boldsymbol{\beta}}_{3}\quad &
\quad\widehat{\boldsymbol{\beta}}_{12}\quad & \quad\widehat{\boldsymbol{\beta}}_{13}\quad &
\quad\widehat{\boldsymbol{\beta}}_{23}\quad
  \end{array}
$$}
\vspace{-0.5cm}
$$
  \widehat{\mathbb{B}} =
  \left [
    \begin{array}{c}
      \widehat{\boldsymbol{\beta}}'_{1} \\
      \widehat{\boldsymbol{\beta}}'_{2}
    \end{array}
  \right ]' =
  \left[
     \begin{array}{ccccccc}
       104.86 & -3.147 & -0.142 & -0.199 & 2.379 & -0.35 & -0.106 \\
       70.45 & -0.348 & 3.59 & 0.28 & 0.323 & -0.45 & 0.614
     \end{array}
  \right]'
$$
Also,
$$
  \widehat{\mathbf{\Sigma}} =
  \left [
    \begin{array}{l r}
      4.190 & 3.546 \\
      3.546 & 4.666
    \end{array}
  \right ]
$$
From where
\begin{eqnarray*}
  \widehat{\Cov}(\Vec \widehat{\mathbb{B}}) &=& \widehat{\mathbf{\Sigma}} \otimes (\mathbf{X}'\mathbf{X})^{-1} \\
&=& \left [
    \begin{array}{cc}
      4.190 & 3.546 \\
      3.546 & 4.666
    \end{array}
  \right ] \otimes 0.03125  \ \mathbf{I}_{7}.
\end{eqnarray*}
In particular,  $\widehat{\Cov}(\widehat{\boldsymbol{\beta}}_{1}) = 0.131 \ \mathbf{I}_{7}$ and
$\widehat{\Cov}(\widehat{\boldsymbol{\beta}}_{2}) = 0.145 \ \mathbf{I}_{7}$. Therefore,  the
estimator of the covariance matrix of response surfaces according to equation (\ref{eq7}) is
$$
  \widehat{\Cov}(\mathbf{Y}(\mathbf{x})) = (1 + x_{1}^{2} + x_{2}^{2} + x_{3}^{2}
  + x_{1}^{2}x_{2}^{2} + x_{1}^{2}x_{3}^{2} + x_{2}^{2}x_{3}^{2})
  \left [
    \begin{array}{cc}
      0.131 & 0.111 \\
      0.111 & 0.145
    \end{array}
    \right ].
$$
Next,  we propose diverse multiobjective stochastic solutions and their deterministic
equivalent corresponding programs:

\begin{itemize}
   \item \emph{Equivalent multiobjective V-model}
    $$
      \build{\min}{}{\mathbf{x} \in \mathfrak{X}}
        \left[
        \begin{array}{cc}
            \widehat{\Var}\left(\widehat{Y}_{1} \left(\mathbf{x}, \widehat{\boldsymbol{\beta}_{1}}\right)\right)
& \widehat{\Cov}\left(\widehat{Y}_{1} \left(\mathbf{x},
\widehat{\boldsymbol{\beta}_{1}}\right),
            \widehat{Y}_{2} \left(\mathbf{x}, \widehat{\boldsymbol{\beta}_{2}}\right)\right) \\
            \widehat{\Cov}\left(\widehat{Y}_{1} \left(\mathbf{x}, \widehat{\boldsymbol{\beta}_{1}}\right),
            \widehat{Y}_{2} \left(\mathbf{x}, \widehat{\boldsymbol{\beta}_{2}}\right)\right)
& \widehat{\Var}\left(\widehat{Y}_{2} \left(\mathbf{x},
\widehat{\boldsymbol{\beta}_{2}}\right)\right)
       \end{array}
     \right].
   $$
   And its corresponding deterministic equivalent program is
   $$
      \build{\min}{}{\mathbf{x} \in \mathfrak{X}}\mathbf{z}'(\mathbf{x})(\mathbf{X}'
      \mathbf{X})^{-1}\mathbf{z}(\mathbf{x}).
   $$
   \item \emph{Equivalent deterministic multiobjective modified E-model}
   $$
     \build{\min}{}{\mathbf{x} \in \mathfrak{X}}
     \left[
     \begin{array}{c}
        \widehat{\mathbf{Y}}\left(\mathbf{x}, \widehat{\mathbb{B}}\right) \\
        \left(\widehat{\Cov}\left(\widehat{\mathbf{Y}}\left(\mathbf{x}, \widehat{\mathbb{B}}\right)\right)\right)^{1/2}
     \end{array}
     \right ]
   $$
   where
   $$
     \widehat{\mathbf{Y}}\left(\mathbf{x}, \widehat{\mathbb{B}}\right)=
        \left(
        \begin{array}{c}
          \widehat{Y}_{1} \left(\mathbf{x}, \widehat{\boldsymbol{\beta}_{1}}\right) \\
          \widehat{Y}_{2} \left(\mathbf{x}, \widehat{\boldsymbol{\beta}_{2}}\right)
        \end{array}
        \right),
   $$
   and
   $\left(\widehat{\Cov}\left(\widehat{\mathbf{Y}}\left(\mathbf{x}, \widehat{\mathbb{B}}\right)\right)\right)^{1/2}$
   is
   $$
        \left(
        \begin{array}{cc}
            \widehat{\Var}\left(\widehat{Y}_{1} \left(\mathbf{x}, \widehat{\boldsymbol{\beta}_{1}}\right)\right)
& \widehat{\Cov}\left(\widehat{Y}_{1} \left(\mathbf{x},
\widehat{\boldsymbol{\beta}_{1}}\right),
            \widehat{Y}_{2} \left(\mathbf{x}, \widehat{\boldsymbol{\beta}_{2}}\right)\right) \\
            \widehat{\Cov}\left(\widehat{Y}_{1} \left(\mathbf{x}, \widehat{\boldsymbol{\beta}_{1}}\right),
            \widehat{Y}_{2} \left(\mathbf{x}, \widehat{\boldsymbol{\beta}_{2}}\right)\right)
& \widehat{\Var}\left(\widehat{Y}_{2} \left(\mathbf{x},
\widehat{\boldsymbol{\beta}_{2}}\right)\right)
       \end{array}
       \right)^{1/2},
   $$
   In this case the deterministic equivalent program can be stated as (among many other options,  including
   lexicographic or $\epsilon$-constraint models)
   $$
     \build{\min}{}{\mathbf{x} \in \mathfrak{X}}
     r_{1}\ f\left(\widehat{\mathbf{Y}}\left(\mathbf{x}, \widehat{\mathbb{B}}\right)\right) +
     r_{2}\ g\left(\left(\widehat{\Cov}\left(\widehat{\mathbf{Y}}\left(\mathbf{x},
     \widehat{\mathbb{B}}\right)\right)\right)^{1/2}\right),
   $$
   where $r_{j} \geq 0$,  $j = 1, 2$ are constants such that $r_{1} + r_{2} = 1$ (in general),
   whose values indicate the relative importance of the expectation and matrix covariance of
   $\widehat{\mathbf{Y}}\left(\mathbf{x}, \widehat{\mathbb{B}}\right)$,  and $f$ and $g$ are value functions.

   In particular,  using (\ref{equi2}),  a deterministic equivalent program via the weighting method is
   $$
     \build{\min}{}{\mathbf{x} \in \mathfrak{X}}
     r_{1}\ \left(w_{1}\ z'(\mathbf{x})\widehat{\boldsymbol{\beta}}_{1} + w_{2}\
     z'(\mathbf{x})\widehat{\boldsymbol{\beta}}_{2} \right) +
     r_{2}\ z'(\mathbf{x})(\mathbf{X}'\mathbf{X})^{-1}z(\mathbf{x}),
   $$
   and assuming that the primary objective function is $g$,  via the $\epsilon$-constraint model we have
   $$
     \begin{array}{c}
       \build{\min}{}{\mathbf{x} \in \mathfrak{X}}
       z'(\mathbf{x})(\mathbf{X}'\mathbf{X})^{-1}z(\mathbf{x})\\
       \mbox{subject to}\\
       z'(\mathbf{x})\widehat{\boldsymbol{\beta}}_{1} = \tau_{1}\\
       z'(\mathbf{x})\widehat{\boldsymbol{\beta}}_{2} = \tau_{2}.\\
     \end{array}
   $$

   \item \emph{Equivalent deterministic multiobjective P-model}
   $$
    \build{\min}{}{\mathbf{x} \in \mathfrak{X}}
    \left(%
    \begin{array}{c}\displaystyle
      \frac{\tau_{1} - z'(\mathbf{x})\widehat{\boldsymbol{\beta}}_{1}}{
        \sqrt{\widehat{\sigma}_{11} \ z'(\mathbf{x})(\mathbf{X}'\mathbf{X})^{-1}z(\mathbf{x})}} \\
      \displaystyle
      \frac{\tau_{2} - z'(\mathbf{x})\widehat{\boldsymbol{\beta}}_{2}}{
        \sqrt{\widehat{\sigma}_{22} \ z'(\mathbf{x})(\mathbf{X}'\mathbf{X})^{-1}z(\mathbf{x})}} \\
    \end{array}%
    \right).
   $$
   In this case the deterministic equivalent program via the weighting method is
   $$
     \build{\min}{}{\mathbf{x} \in \mathfrak{X}}
       w_{1}\left\{\frac{\tau_{1} - z'(\mathbf{x})\widehat{\boldsymbol{\beta}}_{1}}{
        \sqrt{\widehat{\sigma}_{11} \ z'(\mathbf{x})(\mathbf{X}'\mathbf{X})^{-1}z(\mathbf{x})}}\right\} +
       w_{2}\left\{\frac{\tau_{2} - z'(\mathbf{x})\widehat{\boldsymbol{\beta}}_{2}}{
        \sqrt{\widehat{\sigma}_{22} \ z'(\mathbf{x})(\mathbf{X}'\mathbf{X})^{-1}z(\mathbf{x})}}\right\}.
   $$
   And assuming that $\widehat{Y}_{2} \left(\mathbf{x}, \widehat{\boldsymbol{\beta}}_{2}\right)$ is
   the primary objective function,  the deterministic equivalent program via the
   $\epsilon$-constraint method is
   $$
     \begin{array}{c}
    \build{\min}{}{\mathbf{x}\in \mathfrak{X}} \displaystyle\frac{\tau_{2} - z'(\mathbf{x})
    \widehat{\boldsymbol{\beta}}_{2}}{
        \sqrt{\widehat{\sigma}_{22} \ z'(\mathbf{x})(\mathbf{X}'\mathbf{X})^{-1}z(\mathbf{x})}}\\
    \mbox{subject to}\\
      \displaystyle\frac{\tau_{1} - z'(\mathbf{x})\widehat{\boldsymbol{\beta}}_{1}}{
        \sqrt{\widehat{\sigma}_{11} \ z'(\mathbf{x})(\mathbf{X}'\mathbf{X})^{-1}z(\mathbf{x})}} = \tau_{1}.
  \end{array}
   $$
   \item \emph{Equivalent deterministic multiobjective Kataoka model}
   $$
     \build{\min}{}{\mathbf{x} \in \mathfrak{X}}
    \left(%
    \begin{array}{c}
      z'(\mathbf{x}) \widehat{\boldsymbol{\beta}_{1}} + \boldsymbol{\Phi}^{-1}(\delta) \
        \sqrt{\widehat{\sigma}_{11} \ z'(\mathbf{x})(\mathbf{X}'\mathbf{X})^{-1}z(\mathbf{x})}\\
      z'(\mathbf{x}) \widehat{\boldsymbol{\beta}_{2}} + \boldsymbol{\Phi}^{-1}(\delta) \
        \sqrt{\widehat{\sigma}_{22} \ z'(\mathbf{x})(\mathbf{X}'\mathbf{X})^{-1}z(\mathbf{x})}\\
    \end{array}%
    \right).
   $$
   The deterministic equivalent program via the weighting method is
   $$
     \begin{array}{c}
     \build{\min}{}{\mathbf{x} \in \mathfrak{X}}
      w_{1}\left\{z'(\mathbf{x}) \widehat{\boldsymbol{\beta}_{1}} + \boldsymbol{\Phi}^{-1}(\delta) \
        \sqrt{\widehat{\sigma}_{11} \ z'(\mathbf{x})(\mathbf{X}'\mathbf{X})^{-1}z(\mathbf{x})}\right\}\\
      \qquad \qquad + \ w_{2}\left\{z'(\mathbf{x}) \widehat{\boldsymbol{\beta}_{2}} + \boldsymbol{\Phi}^{-1}(\delta) \
        \sqrt{\widehat{\sigma}_{22} \ z'(\mathbf{x})(\mathbf{X}'\mathbf{X})^{-1}z(\mathbf{x})}\right\}\\
     \end{array}
   $$
   Now assuming that $\widehat{Y}_{2} \left(\mathbf{x}, \widehat{\boldsymbol{\beta}}_{2}\right)$ is
   the primary objective function,  the deterministic equivalent program via the
   $\epsilon$-constraint method is
   $$
     \begin{array}{c}
     \build{\min}{}{\mathbf{x}\in \mathfrak{X}} z'(\mathbf{x}) \widehat{\boldsymbol{\beta}_{2}} +
        \boldsymbol{\Phi}^{-1}(\delta) \ \sqrt{\widehat{\sigma}_{22} \ z'(\mathbf{x})(\mathbf{X}'
        \mathbf{X})^{-1}z(\mathbf{x})}\\
     \mbox{subject to}\\
      z'(\mathbf{x}) \widehat{\boldsymbol{\beta}_{1}} + \boldsymbol{\Phi}^{-1}(\delta) \
        \sqrt{\widehat{\sigma}_{11} \ z'(\mathbf{x})(\mathbf{X}'\mathbf{X})^{-1}z(\mathbf{x})} = \tau_{1}.
     \end{array}
   $$
 \end{itemize}

Table \ref{table5} shows the solution of (\ref{equiv1}) by diverse multiobjective stochastic
methods and other methods described in the literature.
\begin{landscape}
\begin{table}[h]
  \centering
  \caption{Comparison of the results of the proposed model and those derived by other methods.}\label{table5}
  \begin{scriptsize}
  \begin{minipage}[t]{460pt}\hspace{-2cm}
  \begin{tabular}{c|crrrrccccc}
    \hline\hline
    \multicolumn{1}{c}{}
& \begin{tabular}{c}
          \textbf{Stochastic} \\
          \textbf{programming} \\
          \textbf{method}
        \end{tabular}
& $x_{1}$ & $x_{2}$ & $x_{3}$ & $F(\mathbf{x})$\footnote{\tiny Objective function} &
$\widehat{Y}_{1}(\mathbf{x})$ & $\widehat{Y}_{2}(\mathbf{x})$ &
$\widehat{\Var}(\widehat{Y}_{1}(\mathbf{x}))$ & $\widehat{\Var}(\widehat{Y}_{2}(\mathbf{x}))$
& $\widehat{\Cov}(\widehat{Y}_{1}(\mathbf{x}), \widehat{Y}_{2}(\mathbf{x}))$\\
    \hline\hline
    \multicolumn{1}{c}{\textbf{\citet{chh:01}}} & -- & 1.000 & 1.000 & -1.000 & -- & 104.612 & 73.574 & 0.917 & 1.021 & 0.776\\
    \hline
& \textbf{Distance Based}
        \footnote{\tiny Squared euclidean distance} & 0.953 & 0.709 & 0.407 & -- & 103.247 & 73.000 & 0.428 & 0.476 & 0.362\\
& \textbf{Robust E-model} & 1.000 & 0.707 & 0.483 & - & 103.332 & 73.000 & 0.470 & 0.523 & 0.397\\
    \textbf{\citet{hbdn:10}\footnote{\tiny Where weights $\mathbf{w}$ are considered deterministic}}& \begin{tabular}{c}
         \textbf{Lexicographic} \\
         (\textbf{First} $\widehat{E}(F(\mathbf{x}))$) \\
       \end{tabular}
& 1.000 & 0.707 & 0.483 & -- & 103.000 & 73.000 & 0.470 & 0.523 & 0.397\\
& \begin{tabular}{c}
         \textbf{Lexicofraphic} \\
         (\textbf{First} $\widehat{\Var}(F(\mathbf{x}))$) \\
       \end{tabular}
& 0.000 & 0.000 & 0.000 & -- & 104.865 & 70.453 & 0.131 & 0.146 & 0.111\\
& \textbf{Modified V-model} & 0.000 & 0.000 & 0.000 & - & 104.865 & 70.453 & 0.131 & 0.146 & 0.111\\
     \hline
& \textbf{V-model}& 0.000 & 0.000 & 0.000 & 1 & 104.865 & 70.453 & 0.131 & 0.146 & 0.111\\
&  \begin{tabular}{c}
         \textbf{Modified E-model}\footnote{\tiny $r_{1} = r_{2} = 0.5$} \\
         (\textbf{Weighting method}) \\
       \end{tabular}
& 0.522 & -1.000 & 0.108 & 39.588 & 102.100 & 66.449 & 0.336 & 0.375 & 0.285\\
     \textbf{Multiobjective}& \begin{tabular}{c}
         \textbf{Modified E-model} \\
         (\textbf{$\epsilon$-constraint}) \\
       \end{tabular}
& 1.000 & 0.707 & 0.452 & 3.511 & 103.019 & 72.992 & 0.460 & 0.512 & 0.389\\
     \textbf{Stochastic}&  \begin{tabular}{c}
         \textbf{P-model} \\
         (\textbf{Weighting method}) \\
       \end{tabular}
& -0.349 & 1.000 & 0.548 & -2.672 & 104.893 & 74.630 & 0.377 & 0.420 & 0.320\\
     \textbf{approaches}& \begin{tabular}{c}
         \textbf{P-model} \\
         (\textbf{$\epsilon$-constraint}) \\
       \end{tabular}
& 0.910 & -0.658 & 0.000 & 8.799 & 100.672 & 67.577 & 0.343 & 0.382 & 0.290\\
& \begin{tabular}{c}
         \textbf{Kataoka}\footnote{\tiny $\delta = 0.95$} \\
         (\textbf{Weighting method}) \\
       \end{tabular}
& 1.000 & -1.000 & 1.000 & 74.989 & 99.039 & 65.405 & 0.917 & 1.021 & 0.776\\
& \begin{tabular}{c}
         \textbf{Kataoka} \\
         (\textbf{$\epsilon$-constraint}) \\
       \end{tabular}
& 0.541 & -1.000 & 0.851 & 67.296 & 101.780 & 66.006 & 0.556 & 0.619 & 0.470\\
& \begin{tabular}{c}
         \textbf{Goal Programming} \\
            \end{tabular}
& 0.844 & 0.605 & 1 & 0 &    102.78 & 72.78 & 0.6764 & 0.6441 & 0.4895\\
    \hline\hline
  \end{tabular}
  \end{minipage}
  \end{scriptsize}
\end{table}
\end{landscape}

Note that while all of these optimisation techniques essentially provide the solution to the
same practical problem,  i.e. that of obtaining the critical value of the variables
$\mathbf{x}$,  from a mathematical point of view and more precisely from the standpoint of
mathematical programming,  problem (\ref{eq8}) as solved by \citet{chh:01} and \citet{hbdn:10}
and problem (\ref{equiv2}) examined in the present paper are not the same. Therefore,  unless a
reasonable basis for comparison is proposed,  Table 2 should be taken simply as an example of
different approaches and different solutions to the practical problem. From the latter,
experts,  researchers or decision makers can select the most suitable method for solving their
own problem in terms of the particular context.

\section*{Conclusions}

It should be emphasised that even a comparison between the diverse techniques proposed in this
paper should be made with appropriate reservations,  since the solutions discussed refer to
different decision-making criteria. For example, Table \ref{table5} shows that the methods
designed to minimise the variance actually obtain a lower variance than the other solutions,
but perhaps the optimum response variables are a little further from the target. Likewise,  the
minimum risk models provide more conservative solutions.

\section*{Acknowledgments}

This work was supported  by the University of Medell\'{\i }n (Medell\'{\i }n,  Colombia) and
Universidad Aut\'onoma Agraria Antonio Narro (M\'{e}xico),   joint grant No. 469,  SUMMA group.
Also,  the first author was partially supported by IDI-Spain,  Grants No. FQM2006-2271 and
MTM2008-05785. This paper was written during J. A. D\'{\i}az-Garc\'{\i}a's stay as a visiting
professor at the Department of Statistics and O. R. of the University of Granada,  Spain.

\end{document}